\documentclass[a4paper,12pt]{article}
\usepackage{pgfplots} 
\pgfplotsset{compat=1.5}
\usepackage{multicol}

\usepackage{hyperref}
\usepackage{placeins}
\usepackage{fullpage}
\usepackage{graphicx}
\usepackage{amsmath}
\usepackage{amssymb}
\usepackage{subfigure}
\usepackage{epsfig}

\usepackage{caption,color}
\usepackage{amsthm}

\usepackage{bm}

\usepackage{url}

\usepackage{accents}

\usepackage{lipsum}

\newcommand\blfootnote[1]{%
  \begingroup
  \renewcommand\thefootnote{}\footnote{#1}%
  \addtocounter{footnote}{-1}%
  \endgroup
}

\theoremstyle{plain}
\newtheorem{thm}{Theorem}[section]
\newtheorem{cor}[thm]{Corollary}
\newtheorem{prop}[thm]{Proposition}

\theoremstyle{remark}

\author{
Chiun-Chuan Chen, Li-Chang Hung$^{\ast}$ \;and Chen-Chih Lai 
\vspace{10mm}
\\
 \small
 Department of Mathematics, National Taiwan University, Taiwan
\\
}

\title{An N-barrier maximum principle for autonomous systems of $n$ species and its application to problems arising from population dynamics
}

\date{
\small
L.-C. Hung dedicates this work to Mach Nguyet Minh
}

\begin{document}
\maketitle



\blfootnote{$^{\ast}$Corresponding author's email address: \texttt{lichang.hung@gmail.com}}


\begin{abstract}

We show that the N-barrier maximum principle (NBMP) remains true for $n$ $(n>2)$ species. In addition, a stronger lower bound in NBMP is given by employing an improved tangent line method. As an application of NBMP, we establish a nonexistence result for traveling wave solutions to the four species Lotka-Volterra system. 










\end{abstract}





\section{Introduction}

The main purpose of the present paper is to establish the N-barrier maximum principle (NBMP) for $n (n>2)$ species, the case $n=2$ having been considered previously (\cite{Chen&Hung16NBMP},\cite{Hung15NBMP}). To be more specific, we study the autonomous system of $n$ species
\begin{equation}\label{eqn: autonomous system of n species}
d_i\,(u_i)_{xx}+\theta\,(u_i)_{x}+u_i^{m_i}\,f_i(u_1,u_2,...,u_n)=0, \quad x\in\mathbb{R}, \quad i=1,2,...,n,
\end{equation}
where $u_i=u_i(x)$, $d_i,m_i>0$, and $f_i(u_1,u_2,...,u_n)\in C^0(\mathbb{R^{+}}\times\mathbb{R^{+}}\times...\times\mathbb{R^{+}})$ for $i=1,2,...,n$; $\theta\in\mathbb{R}$.
Throughout, we assume, unless otherwise stated, that the following hypothesis on $f_i(u_1,u_2,...,u_n)$ is satisfied:
\begin{itemize}
\item [$\mathbf{[H]}$]
For $i=1,2,...,n$, there exist $\bar{u}_i>\underaccent\bar{u}_i>0$  such that
\begin{eqnarray*}
f_i(u_1,u_2,...,u_n)\geq 0 &\text{ whenever } (u_1,u_2,...,u_n)\in \underaccent\bar{\mathcal{R}};\\\\
f_i(u_1,u_2,...,u_n)\leq  0 &\text{ whenever } (u_1,u_2,...,u_n)\in \bar{\mathcal{R}},
\end{eqnarray*}
where 
\begin{eqnarray*}
\underaccent\bar{\mathcal{R}}&=&\Big\{ (u_1,u_2,...,u_n)\;\Big|\; \sum_{i=1}^{n}\frac{\displaystyle u_i}{\displaystyle\underaccent\bar{u}_i}\le 1,\; u_1,u_2,...,u_n\ge 0 \Big\};\\
\bar{\mathcal{R}}&=&\Big\{ (u_1,u_2,...,u_n)\;\Big|\; \sum_{i=1}^{n}\frac{\displaystyle u_i}{\displaystyle\bar{u}_i}\ge 1,\; u_1,u_2,...,u_n\ge 0 \Big\}.
\end{eqnarray*}
\end{itemize}
We couple \eqref{eqn: autonomous system of n species} with the prescribed Dirichlet conditions at $x=\pm\infty$:
\begin{equation}\label{eqn: autonomous system of n species BC}
(u_1,u_2,...,u_n)(-\infty)=\textbf{e}_{-},\quad (u_1,u_2,...,u_n)(\infty)=\textbf{e}_{+},
\end{equation}
where 
\begin{equation}\label{eqn: e- and e+}
\textbf{e}_{-}, \textbf{e}_{+}\in\Big\{ (u_1,u_2,...,u_n) \;\Big|\; u_i^{m_i}\,f_i(u_1,u_2,...,u_n)=0\;  (i=1,2,...,n), u_1,u_2,...,u_n\ge 0\Big\}
\end{equation}
are the equilibria of \eqref{eqn: autonomous system of n species} which connect the solution $(u_1,u_2,...,u_n)(x)$ at $x=-\infty$ and $x=\infty$. This leads to the boundary value problem of \eqref{eqn: autonomous system of n species} and \eqref{eqn: autonomous system of n species BC}:
\begin{equation*}
\textbf{(BVP)}
\begin{cases}
\vspace{3mm}
d_i\,(u_i)_{xx}+\theta\,(u_i)_{x}+u_i^{m_i}\,f_i(u_1,u_2,...,u_n)=0, \quad x\in\mathbb{R}, \quad i=1,2,...,n, \\
(u_1,u_2,...,u_n)(-\infty)=\textbf{e}_{-},\quad (u_1,u_2,...,u_n)(\infty)=\textbf{e}_{+}.
\end{cases}
\end{equation*}
\textbf{(BVP)} arises from the study of traveling waves in the following reaction-diffusion system of $n$ species (\cite{Murray93Mbiology},\cite{Volpert94}):
\begin{equation}\label{eqn: R-D system}
(\omega_i)_t=d_i\,(\omega_i)_{yy}+\omega_i^{m_i}\,f_i(\omega_1,\omega_2,...,\omega_n), \quad y\in\mathbb{R},\;t>0,\quad i=1,2,...,n,
\end{equation}
where $\omega_i(y,t)$ $(i=1,2,...n)$ is the density of the $i$-th species and $d_i$ $(i=1,2,...n)$ represents the diffusion rate of the $i$-th species. A special solution $(u_1(x),u_2(x),...,u_n(x))=(\omega_1(y,t),\omega_2(y,t),...,\omega_n(y,t))$, $x=y-\theta\,t$, where $\theta$ is the propagation speed of the traveling wave, is a traveling wave solution of \eqref{eqn: R-D system}. It is easy to see that the traveling wave solution of such form satisfies \eqref{eqn: autonomous system of n species}. 

Our main result is that \textbf{(BVP)} enjoys the following N-barrier maximum principle.
\vspace{3mm}
   
\begin{thm}[\textbf{NBMP for $n$ Species}]\label{thm: NBMP for n species}
Assume that $\mathbf{[H]}$ holds. Given any set of $\alpha_i>0$ $(i=1,2,...,n)$, suppose that $(u_1(x),u_2(x),...,u_n(x))$ is a nonnegative $C^2$ solution to \textbf{(BVP)}. Then
\begin{equation}\label{eqn: upper and lower bounds of p_generalized}
\underaccent\bar{\lambda}
\leq \sum_{i=1}^{n} \alpha_i\,u_i(x) \leq
\bar{\lambda}, \quad x\in\mathbb{R},
\end{equation}
where
\begin{eqnarray}
\bar{\lambda} & = & 
\Big(
\max_{i=1,2,...,n}
\alpha_i\,\bar{u}_i
\Big)
\Big(
\max_{i=1,2,...,n} d_i 
\Big)
\Big(
\min_{i=1,2,...,n} d_i
\Big)^{-1},
\\
\underaccent\bar{\lambda} & = & 
\Big(
\min_{i=1,2,...,n}
\alpha_i\,\underaccent\bar{u}_i
\Big)
\Big(
\min_{i=1,2,...,n} d_i 
\Big)
\Big(
\max_{i=1,2,...,n} d_i
\Big)^{-1}
\chi,   
\end{eqnarray}
with $\chi$ defined by
\begin{equation}\label{eqn: chi 1 or 0}
\chi
=
\begin{cases}
\vspace{3mm}
0,
\quad \text{if} \quad  \text{\bf e}_{+}=(0,...,0) \quad \text{or} \quad \text{\bf e}_{-}=(0,...,0),\\
1,
\quad \text{otherwise}.
\end{cases}
\end{equation}

\end{thm}

We note that both the lower bound $\underaccent\bar{\lambda}$ and the upper bound $\bar{\lambda}$ in NBMP do not depend explicitly on the propagation speed $\theta$. To illustrate Theorem~\ref{thm: NBMP for n species}, we present an example. For $n=3$, suppose that $m_i=1$ and $f_i(u_1,u_2,u_3)=u_i\,(\sigma_i-c_{i1}\,u_1-c_{i2}\,u_2-c_{i3}\,u_3)$ for $i=1,2,3$. Then \textbf{(BVP)} becomes
\begin{equation}\nonumber
\textbf{(LV3)}
\begin{cases}
\vspace{3mm}
d_1\,(u_1)_{xx}+\theta \,(u_1)_x+u_1\,(\sigma_1-c_{11}\,u_1-c_{12}\,u_2-c_{13}\,u_3)=0, \quad x\in\mathbb{R}, \\
\vspace{3mm}
d_2\,(u_2)_{xx}+\theta \,(u_2)_x+u_2\,(\sigma_2-c_{21}\,u_1-c_{22}\,u_2-c_{23}\,u_3)=0, \quad x\in\mathbb{R},\\
\vspace{3mm}
d_3\,(u_3)_{xx}+\theta \,(u_3)_x+u_3\,(\sigma_3-c_{31}\,u_1-c_{32}\,u_2-c_{33}\,u_3)=0, \quad x\in\mathbb{R},\\
(u_1,u_2,u_3)(-\infty)=\textbf{e}_{-},\quad (u_1,u_2,u_3)(\infty)=\textbf{e}_{+},
\end{cases}
\end{equation}
where
\begin{equation}\label{eqn: e- and e+ 3 species}
\textbf{e}_{-}, \textbf{e}_{+}\in\Big\{ (u_1,u_2,u_3) \;\Big|\; u_i\,(\sigma_i-c_{i1}\,u_1-c_{i2}\,u_2-c_{i3}\,u_3)=0\;  (i=1,2,3), u_1,u_2,u_3\ge 0\Big\}.
\end{equation}
The parameters $d_i$, $\sigma_i$, $c_{ii}$ $(i=1,2,3)$, and $c_{ij}$ $(i,j=1,2,3  \ \text{with} \;i\neq j)$, which are all positive constants, stand for the diffusion rates, intrinsic growth rates, intra-specific competition rates, and inter-specific competition rates, respectively. A solution $(u_1(x),u_2(x),u_3(x))$ to \textbf{(LV3)} is a traveling wave solution which solves the competitive Lotka-Volterra systems of three competing species:
\begin{equation}\label{eqn: L-V systems of three species RD}
(\omega_i)_t=d_i\,(\omega_i)_{yy}+\omega_i\,(\sigma_i-c_{i1}\,\omega_1-c_{i2}\,\omega_2-c_{i3}\,\omega_3), \quad y\in\mathbb{R},\;t>0,\quad i=1,2,3,
\end{equation}
where $(\omega_1(y,t),\omega_2(y,t),\omega_3(y,t))=(u_1(x),u_2(x),u_3(x))$, $x=y-\theta\,t$. When the diffusion terms are absent, \eqref{eqn: L-V systems of three species RD} is the celebrated May-Leonard model (\cite{May&Leonard75}) under the assumption that $\sigma_i=c_{ii}=1$ $(i=1,2,3)$, $c_{12}=c_{23}=c_{31}=\mu_1>0$ and $c_{13}=c_{21}=c_{32}=\mu_2>0$
\begin{equation}\label{eqn: May-Leonard model}
\begin{cases}
\vspace{3mm}
(\omega_1)_t=\omega_1\,(1-\omega_1-\mu_1\,\omega_2-\mu_2\,\omega_3),\quad t>0, \\
\vspace{3mm}
(\omega_2)_t=\omega_2\,(1-\mu_2\,\omega_1-\omega_2-\mu_1\,\omega_3),\quad t>0, \\
(\omega_3)_t=\omega_3\,(1-\mu_1\,\omega_1-\mu_2\,\omega_2-\omega_3),\quad t>0, \\
\end{cases}
\end{equation}
where $(\omega_1,\omega_2,\omega_3)=(\omega_1(t),\omega_2(t),\omega_3(t))$.

From the viewpoint of the study of \textit{competitive exclusion} (\cite{Armstrong80Competitive-exclusion}, \cite{Hsu08Competitive-exclusion}, \cite{Hsu-Smith-Waltman96Competitive-exclusion-coexistence-Competitive}, \cite{Jang13Competitive-exclusion-Leslie-Gower-competition-Allee}, \cite{McGehee77Competitive-exclusion}, \cite{Smith94Competition}) or \textit{competitor-mediated coexistence} (\cite{CantrellWard97Competition-mediatedCoexistence}, \cite{Kastendiek82Competitor-mediatedCoexistence3Species}, \cite{Mimura15DynamicCoexistence3species}), \textbf{(LV3)} or \eqref{eqn: L-V systems of three species RD} arises from investigating problems where one exotic competing species (say, $u_3$) invades the ecological system of two native species (say, $u_1$ and $u_2$) that are competing in the absence of $u_3$. As indicated in \cite{hung2012JJIAM,Kan-on95}, when $u_3(x)$ is absent in \textbf{(LV3)} with $\textbf{e}_{-}=(\frac{\sigma_1}{c_{11}},0)$ and $\textbf{e}_{+}=(0,\frac{\sigma_2}{c_{22}})$, \textbf{(LV3)} under the condition of strong competition (i.e. $\frac{\sigma_1}{c_{11}}>\frac{\sigma_2}{c_{21}}$ and $\frac{\sigma_2}{c_{22}}>\frac{\sigma_1}{c_{12}}$) admits solutions $(u_1(x),u_2(x))$ having profiles with $u_1(x)$ being monotonically decreasing and $u_2(x)$ being monotonically increasing. Since $u_1(x)$ and $u_2(x)$ dominate the neighborhoods of $x=-\infty$ and $x=\infty$, respectively, we are led to expect that the profile of $u_3(x)$ must be \textit{pulse-like} (we call $u_3(x)$ a pulse if $u_3(-\infty)=u_3(\infty)=0$ and $u_3(x)>0$ for $x\in\mathbb{R}$) if it exists since $u_3$ will prevail only when $u_1$ and $u_2$ are not dominant. It turns out that this conjecture is true under certain assumptions on the parameters. In \cite{CHMU-semi,CHMU}, we established existence of this type of solution by finding exact traveling wave solutions in addition to numerical experiments. 

To the best of our knowledge, however, a priori estimates for the parameter dependence of solutions to \textbf{(LV3)} have not yet been found. Corollary~\ref{cor: NBMP for 3 species} provides an affirmative answer to the following question:  

\textbf{Q}: \textit{Can upper and lower bounds of $u_1+u_2+u_3$ can be given in terms of the parameters in \textbf{(LV3)}?}

The above question arises in attempts to understand the ecological capacity of the inhabitant of the three competing species $u_1$, $u_2$, and $u_3$. Due to limited resources, the investigation of the total density of the three species is of interest. More generally, estimates of $\alpha_1\,u_1+\alpha_2\,u_2+\alpha_3\,u_3$, where $\alpha_i>0$ $(i=1,2,3)$, are given in Corollary~\ref{cor: NBMP for 3 species}.

\begin{cor}[\textbf{NBMP for Lotka-Volterra systems of three competing species}]\label{cor: NBMP for 3 species}
Assume that $(u(x),v(x),w(x))$ is a nonnegative $C^2$ solution to \textbf{(LV3)}. For any set of $\alpha_i>0$ $(i=1,2,3)$, we have
\begin{equation}\label{eqn: upper and lower bounds of p_generalized}
\underaccent\bar{\lambda}
\leq \alpha_1\,u_1(x)+\alpha_2\,u_2(x)+\alpha_3\,u_3(x) \leq
\bar{\lambda}, \quad x\in\mathbb{R},
\end{equation}
where
\begin{eqnarray}
\bar{\lambda} & = & \max_{i=1,2,3}
\Big(
\alpha_i\max_{j=1,2,3}\frac{\sigma_j}{c_{ji}}
\Big)\,
\Big(
\max_{i=1,2,...,n} d_i 
\Big)
\Big(
\min_{i=1,2,...,n} d_i
\Big)^{-1}, \\
\underaccent\bar{\lambda} & = & \min_{i=1,2,3}
\Big(
\alpha_i\min_{j=1,2,3}\frac{\sigma_j}{c_{ji}}
\Big)\,
\Big(
\min_{i=1,2,...,n} d_i 
\Big)
\Big(
\max_{i=1,2,...,n} d_i
\Big)^{-1}
\,\chi,
\end{eqnarray}
with $\chi$ defined by
\begin{equation}
\chi
=
\begin{cases}
\vspace{3mm}
0,
\quad \text{if} \quad  \text{\bf e}_{+}=(0,0,0) \quad \text{or} \quad \text{\bf e}_{-}=(0,0,0),\\
1,
\quad \text{otherwise}.
\end{cases}
\end{equation}
\end{cor}

\begin{proof}
We apply Theorem~\ref{thm: NBMP for n species} to prove Corollary~\ref{cor: NBMP for 3 species}. Taking
\begin{eqnarray}
\bar{u}_i & = & \max_{j=1,2,3}\frac{\sigma_j}{c_{ji}}; \\
\underaccent\bar{u}_i & = & \min_{j=1,2,3}\frac{\sigma_j}{c_{ji}}.
\end{eqnarray}
It can be verified that $\mathbf{[H]}$ is satisfied. Indeed, we have
\begin{eqnarray*}
\underaccent\bar{\mathcal{R}}&=&\Bigg\{ (u_1,u_2,...,u_n)\;\Bigg|\; \sum_{i=1}^{n}\frac{\displaystyle u_i}{\displaystyle\min_{j=1,2,3}\frac{\sigma_j}{c_{ji}}}\le 1,\; u_1,u_2,...,u_n\ge 0 \Bigg\};\\
\bar{\mathcal{R}}&=&\Bigg\{ (u_1,u_2,...,u_n)\;\Bigg|\; \sum_{i=1}^{n}\frac{\displaystyle u_i}{\displaystyle\max_{j=1,2,3}\frac{\sigma_j}{c_{ji}}}\ge 1,\; u_1,u_2,...,u_n\ge 0 \Bigg\}.
\end{eqnarray*}
Since $\displaystyle \min_{j=1,2,3}\frac{\sigma_j}{c_{ji}}$ ($\displaystyle \max_{j=1,2,3}\frac{\sigma_j}{c_{ji}}$, respectively) is the smallest (largest, respectively) $u_i$-intercept of the three hyperplanes  $\sigma_i-c_{i1}\,u_1-c_{i2}\,u_2-c_{i3}\,u_3=0$ $(i=1,2,3)$, we see that 
\begin{eqnarray*}
\sigma_i-c_{i1}\,u_1-c_{i2}\,u_2-c_{i3}\,u_3\ge 0 &\text{ whenever } (u_1,u_2,...,u_n)\in \underaccent\bar{\mathcal{R}};\\\\
\sigma_i-c_{i1}\,u_1-c_{i2}\,u_2-c_{i3}\,u_3\le 0 &\text{ whenever } (u_1,u_2,...,u_n)\in \bar{\mathcal{R}},
\end{eqnarray*}
for each $i=1,2,3$. The desired result follows from Theorem~\ref{thm: NBMP for n species}.

\end{proof}


NBMP for the diffusive Lotka-Volterra system of two competing species was established in \cite{Chen&Hung16NBMP}, where it was also shown that under additional restrictions on the parameters, a lower bound stronger than the one given in Proposition~\ref{prop: lower bed in NBMP for 2 species} can be found by employing the tangent line method. 

\begin{prop}[\cite{Chen&Hung16NBMP}]
\label{prop: lower bed in NBMP for 2 species}
Let $a_1>1$ and $a_2>1$. Suppose that $(u(x),v(x))$ is $C^2$, nonnegative, and satisfies the following differential inequalities
and asymptotic behavior:
\begin{equation}\label{eqn: L-V <0}
\begin{cases}
\vspace{3mm}
\hspace{2.0mm} u_{xx}+\theta\,u_{x}+u\,(1-u-a_1\,v)\leq0, \quad x\in\mathbb{R}, \\
\vspace{3mm}
d\,v_{xx}+\theta\,v_{x}+k\,v\,(1-a_2\,u-v)\leq0, \quad x\in\mathbb{R},\\
(u,v)(-\infty)=(1,0),\quad (u,v)(+\infty)=(0,1),
\end{cases}
\end{equation}
where $d$, $k$, $a_1$, $a_2$ are positive constants. For any $\alpha, \beta>0$, we have 
\begin{equation}\label{q lower bound}
\alpha\,u(x)+\beta\,v(x)\geq \min\bigg[\frac{\alpha}{a_2\,d},\frac{\beta}{a_1}\bigg]\,\min [1,d^2],\quad x\in\mathbb{R}.
\end{equation}

\end{prop}

We show in Section~\ref{sec: stronger lower bound} that the tangent line method can be improved so that the additional parameter restrictions for giving a stronger lower bound than the one given in Proposition~\ref{prop: lower bed in NBMP for 2 species} are no longer needed and, additionally, this lower bound holds for $a_1 >1$ and $a_2 >1$.


The remainder of the paper is organized as follows. NBMP for $n$ species (Theorem~\ref{thm: NBMP for n species}) is proved in Section~\ref{sec: proof of NBMP for n species}. As an application of Corollary~\ref{cor: NBMP for 3 species}, we establish in Section~\ref{sec: nonexistence} a nonexistence result for traveling wave solutions of the Lotka-Volterra system for four competing species
\begin{equation}\nonumber
\textbf{(LV4)}
\begin{cases}
\vspace{3mm}
d_i(u_i)_{xx}+\theta(u_i)_{x}+u_i(\sigma_i-c_{i1}\,u_1-c_{i2}\,u_2-c_{i3}\,u_3-c_{i4}\,u_4)=0,\; x\in\mathbb{R},\; i=1,...,4, \\
(u_1,u_2,u_3,u_4)(-\infty)=(\frac{\sigma_1}{c_{11}},0,0,0),\quad 
(u_1,u_2,u_3,u_4)(\infty)=(0,\frac{\sigma_2}{c_{22}},0,0),
\end{cases}
\end{equation}
where $d_i$, $\sigma_i$, and $c_{ij}$ $(i,j=1,2,3,4)$ are positive constants; $\theta\in\mathbb{R}$ is the propagation speed of the traveling wave.

\vspace{2mm}
\setcounter{equation}{0}
\setcounter{figure}{0}
\setcounter{subfigure}{0}
\section{Proof of Theorem~\ref{thm: NBMP for n species}}\label{sec: proof of NBMP for n species}
\vspace{2mm}

In this section, we prove Theorem~\ref{thm: NBMP for n species}. To this end, we first show in Proposition~\ref{prop: lower bed} that the lower bound given in Theorem~\ref{thm: NBMP for n species} holds when $(u_1,u_2,...,u_n)(x)$ is an upper solution of \textbf{(BVP)} by constructing an appropriate N-barrier.


\vspace{2mm}

\begin{prop} [\textbf{Lower bound in NBMP}]\label{prop: lower bed}
Suppose that $u_i(x)\in C^2(\mathbb{R})$ with $u_i(x)\ge0$ $(i=1,2,...,n)$ and satisfy the following differential inequalities and asymptotic behavior:
\begin{equation*}
\textbf{(Upper)}
\begin{cases}
\vspace{3mm}
d_i\,(u_i)_{xx}+\theta\,(u_i)_{x}+u_i^{m_i}\,f_i(u_1,u_2,...,u_n)\le0, \quad i=1,2,...,n, \quad x\in\mathbb{R},\\
(u_1,u_2,...,u_n)(-\infty)=\textbf{e}_{-},\quad (u_1,u_2,...,u_n)(\infty)=\textbf{e}_{+},
\end{cases}
\end{equation*}
where $\textbf{e}_{-}$ and $\textbf{e}_{+}$ are given by \eqref{eqn: e- and e+}. If the hypothesis 
\begin{itemize}
\item [$\mathbf{[\underline{H}]}$]
For $i=1,2,...,n$, there exist $\underaccent\bar{u}_i>0$  such that
\begin{eqnarray*}
f_i(u_1,u_2,...,u_n)\geq 0 &\text{ whenever } (u_1,u_2,...,u_n)\in \underaccent\bar{\mathcal{R}},
\end{eqnarray*}
where $\underaccent\bar{\mathcal{R}}$ is as defined in $\mathbf{[H]}$
\end{itemize}
holds, then we have for any $\alpha_i>0$ $(i=1,2,...,n)$
\begin{equation}\label{eqn: lower bound of p}
\sum_{i=1}^{n} \alpha_i\,u_i(x)\geq \min\big(\alpha_1\,\underaccent\bar{u}_1,\alpha_2\,\underaccent\bar{u}_2,...,\alpha_n\,\underaccent\bar{u}_n\big)\,\frac{\min(d_1,d_2,...,d_n)}{\max(d_1,d_2,...,d_n)}\,\chi, \quad x\in\mathbb{R},
\end{equation}
where $\chi$ is defined as in \eqref{eqn: chi 1 or 0}.
\end{prop}

\begin{proof}
For the case where $\text{\bf e}_{+}=(0,...,0)$ or $\text{\bf e}_{-}=(0,...,0)$, a trivial lower bound $0$ of $\sum_{i=1}^{n} \alpha_i\,u_i(x)$ is obvious. It suffices to show \eqref{eqn: lower bound of p} for the case $\text{\bf e}_{+}\neq(0,...,0)$ and $\text{\bf e}_{-}\neq(0,...,0)$. To this end, we let
\begin{eqnarray}
p(x) & = & \sum_{i=1}^{n} \alpha_i\,u_i(x);\\
q(x) & = & \sum_{i=1}^{n} \alpha_i\,\,d_i\,u_i(x).
\end{eqnarray}
Adding the $n$ equations in \textbf{(Upper)}, we obtain a single equation involving $p(x)$ and $q(x)$
\begin{equation}\label{eqn: ODE for p and q}
\frac{d^2 q(x)}{dx^2}+\theta\,\frac{d p(x)}{dx}+F(u_1(x),u_2(x),...,u_n(x))\le0,\quad x\in\mathbb{R},
\end{equation}
where $F(u_1,u_2,...,u_n)=\sum_{i=1}^{n} \alpha_i\,u_i^{m_i}\,f_i(u_1,u_2,...,u_n)$. First of all, we treat the case of $d_i\neq d_j$ at least for some $i,j\in \{ 1,2,...,n \}$. 

Determining an appropriate \textit{N-barrier} is crucial in establishing \eqref{eqn: lower bound of p}. The construction of the N-barrier consists of determining $\lambda_2$, $\eta$, and $\lambda_1$ such that the three hyperplanes $\sum_{i=1}^{n} \alpha_i\,\,d_i\,u_i=\lambda_2$, $\sum_{i=1}^{n} \alpha_i\,u_i=\eta$ and $\sum_{i=1}^{n} \alpha_i\,\,d_i\,u_i=\lambda_1$ satisfy the relationship
\begin{equation}\label{eqn: Q1<P<Q2<R}
\mathcal{Q}_{\lambda_1}\subset \mathcal{P}_{\eta}\subset \mathcal{Q}_{\lambda_2}\subset \underaccent\bar{\mathcal{R}},
\end{equation}
where 
\begin{eqnarray}
\mathcal{P}_{\eta} & = & \Big\{ (u_1,u_2,...,u_n) \;\Big|\; \sum_{i=1}^{n} \alpha_i\,u_i\le\eta,\; u_1,u_2,...,u_n\ge 0\Big\}; \\
\mathcal{Q}_{\lambda} & = & \Big\{ (u_1,u_2,...,u_n) \;\Big|\; \sum_{i=1}^{n} \alpha_i\,\,d_i\,u_i\le\lambda,\; u_1,u_2,...,u_n\ge 0\Big\}.
\end{eqnarray}
We follow the three steps below to construct the N-barrier:
\begin{enumerate}
  \item Taking $\lambda_2=\min\{\alpha_1\, d_1\,\underaccent\bar{u}_1,\alpha_2\, d_2\,\underaccent\bar{u}_2,...\alpha_n\, d_n\,\underaccent\bar{u}_n\}$, the hyperplane $\sum_{i=1}^{n} \alpha_i\,\,d_i\,u_i=\lambda_2$ has the $n$ intercepts $(\frac{\displaystyle\lambda_2}{\displaystyle\alpha_1\, d_1},0,...,0)$, $(0,\frac{\displaystyle\lambda_2}{\displaystyle\alpha_2\, d_2},0,...,0)$,..., and $(0,0,...,0,\frac{\displaystyle\lambda_2}{\displaystyle\alpha_n\, d_n})$. It is readily seen that $\frac{\displaystyle\lambda_2}{\displaystyle\alpha_i\, d_i}\le \underaccent\bar{u}_i$ for $i=1,2,...,n$, which gives $\mathcal{Q}_{\lambda_2}\subset \underaccent\bar{\mathcal{R}}$;
  \item Taking $\eta=\lambda_2\,\min\Big\{ \frac{\displaystyle1}{\displaystyle d_1},\frac{\displaystyle1}{\displaystyle d_2},...,\frac{\displaystyle1}{\displaystyle d_n} \Big\}$, the hyperplane $\sum_{i=1}^{n} \alpha_i\,u_i=\eta$ has the $n$ intercepts $(\frac{\displaystyle\eta}{\displaystyle\alpha_1},0,...,0)$, $(0,\frac{\displaystyle\eta}{\displaystyle\alpha_2},0,...,0)$,..., and $(0,0,...,0,\frac{\displaystyle\eta}{\displaystyle\alpha_n})$. It is readily seen that $\frac{\displaystyle\eta}{\displaystyle\alpha_i}\le \frac{\displaystyle\lambda_2}{\displaystyle\alpha_i\, d_i}$ for $i=1,2,...,n$, which gives $\mathcal{P}_{\eta}\subset \mathcal{Q}_{\lambda_2}$;
  \item Taking $\lambda_1=\eta\,\min\{d_1,d_2,...,d_n\}$, the hyperplane $\sum_{i=1}^{n} \alpha_i\,\,d_i\,u_i=\lambda_1$ has the $n$ intercepts $(\frac{\displaystyle\lambda_1}{\displaystyle\alpha_1\, d_1},0,...,0)$, $(0,\frac{\displaystyle\lambda_1}{\displaystyle\alpha_2\, d_2},0,...,0)$,..., and $(0,0,...,0,\frac{\displaystyle\lambda_1}{\displaystyle\alpha_n\, d_n})$. It is readily seen that $\frac{\displaystyle\lambda_1}{\displaystyle\alpha_i\, d_i}\le \frac{\displaystyle\eta}{\displaystyle\alpha_i}$ for $i=1,2,...,n$, which gives $\mathcal{Q}_{\lambda_1}\subset \mathcal{P}_{\eta}$.
\end{enumerate}
The three hyperplanes $\sum_{i=1}^{n} \alpha_i\,\,d_i\,u_i=\lambda_2$, $\sum_{i=1}^{n} \alpha_i\,u_i=\eta$ and $\sum_{i=1}^{n} \alpha_i\,\,d_i\,u_i=\lambda_1$ constructed above form the N-barrier. From the above three steps, it follows immediately that $\lambda_1$ is given by
\begin{equation}\label{eqn: lambda1}
\lambda_1=
\min\big(\alpha_1\,d_1\,\underaccent\bar{u}_1,\alpha_2\,d_2\,\underaccent\bar{u}_2,...,\alpha_n\,d_n\,\underaccent\bar{u}_n\big)\,\frac{\min(d_1,d_2,...,d_n)}{\max(d_1,d_2,...,d_n)}.
\end{equation}
We claim that $q(x)\ge \lambda_1$, $x\in\mathbb{R}$. This proves \eqref{eqn: lower bound of p} since the $\alpha_i>0$ $(i=1,2,...,n)$ are arbitrary. Suppose that, contrary to our claim, there exists $z\in\mathbb{R}$ such that $q(z)<\lambda_1$. Since $u,v\in C^2(\mathbb{R})$ and  $(u_1,u_2,...,u_n)(\pm\infty)=\textbf{e}_{\pm}$, we may assume $\min_{x\in\mathbb{R}} q(x)=q(z)$. 
We denote respectively by $z_2$ and $z_1$ the first points at which the solution $(u_1(x),u_2(x),...,u_n(x))$ intersects the hyperplane $\sum_{i=1}^{n} \alpha_i\,\,d_i\,u_i=\lambda_2$ when $x$ moves from $z$ towards $\infty$ and $-\infty$. 
For the case where $\theta\leq0$, we integrate \eqref{eqn: ODE for p and q} with respect to $x$ from $z_1$ to $z$ and obtain
\begin{equation}\label{eqn: integrating eqn}
q'(z)-q'(z_1)+\theta\,(p(z)-p(z_1))+\int_{z_1}^{z}F(u_1(x),u_2(x),...,u_n(x))\,dx\leq0.
\end{equation}
On the other hand we have:
\begin{itemize}
  \item due to $\min_{x\in\mathbb{R}} q(x)=q(z)$, $q'(z)=0$;
  \item $q(z_1)=\lambda_2$ follows from the fact that $z_1$ is on the hyperplane $\sum_{i=1}^{n} \alpha_i\,\,d_i\,u_i=\lambda_2$. Since $z_1$ is the first point
  for $q(x)$ taking the value $\lambda_2$ when $x$ moves from $z$ to $-\infty$, we conclude that $q(z_1+\delta)\leq \lambda_2$ for $z-z_1>\delta>0$
  and $q'(z_1)\leq 0$;
  \item $p(z)<\eta$ since $z$ is below the hyperplane $\sum_{i=1}^{n} \alpha_i\,u_i=\eta$; $p(z_1)>\eta$ since $z_1$ is above the hyperplane $\sum_{i=1}^{n} \alpha_i\,u_i=\eta$;
  \item let $\mathcal{F}_+=\{(u_1,u_2,...,u_n)\,|\, F(u_1,u_2,...,u_n)> 0, u_1,u_2,...,u_n\ge0\}$. 
  Due to the fact that $(u_1(z_1),u_2(z_1),...,u_n(z_1))$ is on the hyperplane $\sum_{i=1}^{n} \alpha_i\,\,d_i\,u_i=\lambda_2$ and $(u_1(z),u_2(z),...,u_n(z))$ $\in\mathcal{Q}_{\lambda_1}$, $(u_1(z_1),u_2(z_1),...,u_n(z_1))$, $(u_1(z),u_2(z),...,u_n(z))$ $\in\underaccent\bar{\mathcal{R}}$ by \eqref{eqn: Q1<P<Q2<R}. Because of $\mathbf{[\underline{H}]}$ and $F(u_1,u_2,...,u_n)=\sum_{i=1}^{n} \alpha_i\,u_i^{m_i}\,f_i(u_1,u_2,...,u_n)$, it is easy to see that $\{(u_1(x),u_2(x),...,u_n(x))\,|\,z_1\le x \le z\}$ $\subset\underaccent\bar{\mathcal{R}}$ $\subset \mathcal{F}_+$.
  Therefore we have $\int_{z_1}^{z}F(u_1(x),u_2(x),...,u_n(x))\,dx>0$.
\end{itemize}

Combining the above arguments, we obtain
\begin{equation}
q'(z)-q'(z_1)+\theta\,(p(z)-p(z_1))+\int_{z_1}^{z}F(u_1(x),u_2(x),...,u_n(x))\,dx>0,
\end{equation}
which contradicts \eqref{eqn: integrating eqn}. Therefore when $\theta\leq0$, $q(x)\geq \lambda_1$ for $x\in \mathbb{R}$. For the case where $\theta\geq0$, integrating \eqref{eqn: ODE for p and q} with respect to $x$ from $z$ to $z_2$ yields
\begin{equation}\label{eqn: eqn by integrate from z to z2}
q'(z_2)-q'(z)+\theta\,(p(z_2)-p(z))+\int_{z}^{z_2}F(u_1(x),u_2(x),...,u_n(x))\,dx\leq0.
\end{equation}
In a similar manner, it can be shown that $q'(z_2)\ge 0$, $q'(z)=0$, $p(z_2)>\eta$, $p(z)<\eta$, and $\int_{z}^{z_2}F(u_1(x),u_2(x),...,u_n(x))\,dx>0$. These together contradict \eqref{eqn: eqn by integrate from z to z2}. Consequently, \eqref{eqn: lower bound of p} is proved for the case of $d_i\neq d_j$ at least for some $i,j\in \{ 1,2,...,n \}$. Now we turn to the case of $d_i=d$ for all $i=1,2,...,n$. In this case, $q(x)=d\,p(x)$ and \eqref{eqn: ODE for p and q} becomes
\begin{equation}\label{eqn: ODE for p and q with equal di}
d\,\frac{d^2 p(x)}{dx^2}+\theta\,\frac{d p(x)}{dx}+F(u_1(x),u_2(x),...,u_n(x))\le0,\quad x\in\mathbb{R}.
\end{equation}
We take $\lambda_1=\lambda_2=d\,\min\{\alpha_1\,\underaccent\bar{u}_1,\alpha_2\,\underaccent\bar{u}_2,...,\alpha_n\,\underaccent\bar{u}_n\}$ and $\eta=\min\{\alpha_1\,\underaccent\bar{u}_1,\alpha_2\,\underaccent\bar{u}_2,...,\alpha_n\,\underaccent\bar{u}_n\}$. It follows that the three hyperplanes $\sum_{i=1}^{n} \alpha_i\,\,d_i\,u_i=\lambda_2$, $\sum_{i=1}^{n} \alpha_i\,u_i=\eta$ and $\sum_{i=1}^{n} \alpha_i\,\,d_i\,u_i=\lambda_1$ coincide. Analogously to the previous case, we assume that there exists $\hat{z}\in\mathbb{R}$ such that $p(\hat{z})<\lambda_1$ and $\min_{x\in\mathbb{R}} p(x)=p(\hat{z})$. Due to $\min_{x\in\mathbb{R}} p(x)=p(\hat{z})$, we have $p'(\hat{z})=0$ and $p''(\hat{z})\geq0$. Since $(u_1(\hat z),u_2(\hat z),...,u_n(\hat z))$ is in the interior of $\underaccent\bar{\mathcal{R}}$, which is contained in the interior of $\mathcal{F}_+$, we have $F(u_1(\hat z),u_2(\hat z),...,u_n(\hat z))>0$. These together give $d\,p''(\hat{z})+\theta\,p'(\hat{z})+F(u_1(\hat z),u_2(\hat z),...,u_n(\hat z))>0$, which contradicts \eqref{eqn: ODE for p and q with equal di}. As a result, $p(x)\geq \lambda_1$ for $x\in \mathbb{R}$ when $d_i=d$ for all $i=1,2,...,n$. The proof is completed.
\end{proof}
\vspace{3mm}


When $(u_1,u_2,...,u_n)(x)$ is a lower solution of \textbf{(BVP)}, the upper bound in Theorem~\ref{thm: NBMP for n species} can be proved in a similar manner.


\begin{prop} [\textbf{Upper bound in NBMP}]\label{prop: upper bed}
Suppose that $u_i(x)\in C^2(\mathbb{R})$ with $u_i(x)\ge0$ $(i=1,2,...,n)$ and satisfy the following differential inequalities and asymptotic behavior:
\begin{equation*}
\textbf{(Lower)}
\begin{cases}
\vspace{3mm}
d_i\,(u_i)_{xx}+\theta\,(u_i)_{x}+u_i^{m_i}\,f_i(u_1,u_2,...,u_n)\ge0, \quad i=1,2,...,n, \quad x\in\mathbb{R},\\
(u_1,u_2,...,u_n)(-\infty)=\textbf{e}_{-},\quad (u_1,u_2,...,u_n)(\infty)=\textbf{e}_{+},
\end{cases}
\end{equation*}
where $\textbf{e}_{-}$ and $\textbf{e}_{+}$ are given by \eqref{eqn: e- and e+}. If the hypothesis 
\begin{itemize}
\item [$\mathbf{[\bar{H}]}$]
For $i=1,2,...,n$, there exist $\bar{u}_i>0$  such that
\begin{eqnarray*}
f_i(u_1,u_2,...,u_n)< 0 &\text{ whenever } (u_1,u_2,...,u_n)\in\bar{\mathcal{R}},
\end{eqnarray*}
where $\bar{\mathcal{R}}$ is as defined in $\mathbf{[H]}$
\end{itemize}
holds, then we have for any $\alpha_i>0$ $(i=1,2,...,n)$
\begin{equation}\label{eqn: upper bound of p}
\sum_{i=1}^{n} \alpha_i\,u_i(x)\leq \max\big(\alpha_1\,\bar{u}_1,\alpha_2\,\bar{u}_2,...,\alpha_n\,\bar{u}_n\big)\,\frac{\max(d_1,d_2,...,d_n)}{\min(d_1,d_2,...,d_n)}.
\end{equation}
\end{prop}

\begin{proof}
The proof lies in the fact that an appropriate N-barrier for the upper bound \eqref{eqn: upper bound of p} can be constructed. Let
\begin{eqnarray}
\mathcal{P}_{\eta} & = & \Big\{ (u_1,u_2,...,u_n) \;\Big|\; \sum_{i=1}^{n} \alpha_i\,u_i\ge\eta,\; u_1,u_2,...,u_n\ge 0\Big\}; \\
\mathcal{Q}_{\lambda} & = & \Big\{ (u_1,u_2,...,u_n) \;\Big|\; \sum_{i=1}^{n} \alpha_i\,\,d_i\,u_i\ge\lambda,\; u_1,u_2,...,u_n\ge 0\Big\}.
\end{eqnarray}
We determine $\lambda_2$, $\eta$, and $\lambda_1$ in the following steps:
\begin{enumerate}
  \item Taking $\lambda_2=\max\{\alpha_1\, d_1\,\bar{u}_1,\alpha_2\, d_2\,\bar{u}_2,...\alpha_n\, d_n\,\bar{u}_n\}$, the hyperplane $\sum_{i=1}^{n} \alpha_i\,\,d_i\,u_i=\lambda_2$ has the $n$ intercepts $(\frac{\displaystyle\lambda_2}{\displaystyle\alpha_1\, d_1},0,...,0)$, $(0,\frac{\displaystyle\lambda_2}{\displaystyle\alpha_2\, d_2},0,...,0)$,..., and $(0,0,...,0,\frac{\displaystyle\lambda_2}{\displaystyle\alpha_n\, d_n})$. It follows that $\frac{\displaystyle\lambda_2}{\displaystyle\alpha_i\, d_i}\ge \bar{u}_i$ for $i=1,2,...,n$, which gives $\mathcal{Q}_{\lambda_2}\supset \underaccent\bar{\mathcal{R}}$;
  \item Taking $\eta=\lambda_2\,\max\Big\{ \frac{\displaystyle1}{\displaystyle d_1},\frac{\displaystyle1}{\displaystyle d_2},...,\frac{\displaystyle1}{\displaystyle d_n} \Big\}$, the hyperplane $\sum_{i=1}^{n} \alpha_i\,u_i=\eta$ has the $n$ intercepts $(\frac{\displaystyle\eta}{\displaystyle\alpha_1},0,...,0)$, $(0,\frac{\displaystyle\eta}{\displaystyle\alpha_2},0,...,0)$,..., and $(0,0,...,0,\frac{\displaystyle\eta}{\displaystyle\alpha_n})$. It follows that $\frac{\displaystyle\eta}{\displaystyle\alpha_i}\ge \frac{\displaystyle\lambda_2}{\displaystyle\alpha_i\, d_i}$ for $i=1,2,...,n$, which gives $\mathcal{P}_{\eta}\supset \mathcal{Q}_{\lambda_2}$;
  \item Taking $\lambda_1=\eta\,\max\{d_1,d_2,...,d_n\}$, the hyperplane $\sum_{i=1}^{n} \alpha_i\,\,d_i\,u_i=\lambda_1$ has the $n$ intercepts $(\frac{\displaystyle\lambda_1}{\displaystyle\alpha_1\, d_1},0,...,0)$, $(0,\frac{\displaystyle\lambda_1}{\displaystyle\alpha_2\, d_2},0,...,0)$,..., and $(0,0,...,0,\frac{\displaystyle\lambda_1}{\displaystyle\alpha_n\, d_n})$. It follows that $\frac{\displaystyle\lambda_1}{\displaystyle\alpha_i\, d_i}\ge \frac{\displaystyle\eta}{\displaystyle\alpha_i}$ for $i=1,2,...,n$, which gives $\mathcal{Q}_{\lambda_1}\supset \mathcal{P}_{\eta}$.
\end{enumerate}
The three hyperplanes $\sum_{i=1}^{n} \alpha_i\,\,d_i\,u_i=\lambda_2$, $\sum_{i=1}^{n} \alpha_i\,u_i=\eta$ and $\sum_{i=1}^{n} \alpha_i\,\,d_i\,u_i=\lambda_1$ constructed above form the N-barrier which satisfies the property
\begin{equation}\label{eqn: Q1>P>Q2>R}
\mathcal{Q}_{\lambda_1}\supset \mathcal{P}_{\eta}\supset \mathcal{Q}_{\lambda_2}\supset \bar{\mathcal{R}}.
\end{equation}
It follows immediately that $\lambda_1$ is given by
\begin{equation}\label{eqn: lambda1}
\lambda_1=
\max\big(\alpha_1\,\bar{u}_1,\alpha_2\,\bar{u}_2,...,\alpha_n\,\bar{u}_n\big)\,\frac{\max(d_1,d_2,...,d_n)}{\min(d_1,d_2,...,d_n)}.
\end{equation}
We claim that $q(x)=\sum_{i=1}^{n} \alpha_i\,\,d_i\,u_i(x)\le \lambda_1$, $x\in\mathbb{R}$ by contradiction as we have done in Propositions~\ref{prop: lower bed}. The detailed proof of the claim is omitted here for brevity. This completes the proof.
\end{proof}

We are now in the position to prove Theorem~\ref{thm: NBMP for n species}.

\begin{proof}[Proof of Theorem~\ref{thm: NBMP for n species}]
In Propositions~\ref{prop: lower bed} and \ref{prop: upper bed}, we obtain a lower and upper bound for $\sum_{i=1}^{n} \alpha_i\,u_i(x)$, respectively. Combining the results in Propositions~\ref{prop: lower bed} and \ref{prop: upper bed}, we immediately establish Theorem~\ref{thm: NBMP for n species}.
\end{proof}

\vspace{2mm}
\setcounter{equation}{0}
\setcounter{figure}{0}
\setcounter{subfigure}{0}
\section{Improved tangent line method
}\label{sec: stronger lower bound}
\vspace{2mm}

In \cite{Chen&Hung16NBMP}, it is shown that under certain restrictions on the parameters, the lower bound in Proposition~\ref{prop: lower bed in NBMP for 2 species} can be improved by means of the tangent line method. In this section, we show that an improved lower bound can be given without additional conditions on the parameters. To achieve this, let us denote by $\mathcal{L}$ the quadratic curve $\alpha\,u\,(1-u-a_1\,v)+\beta\,k\,v\,(1-a_2\,u-v)=0$ in the first quadrant of the $uv$-plane, i.e.
\begin{equation}
\mathcal{L}=\Big\{(u,v)\,\Big|\, H(u,v)=0, u\ge0, v\ge0\Big\},
\end{equation}
where $ H(u,v):=\alpha\,u\,(1-u-a_1\,v)+\beta\,k\,v\,(1-a_2\,u-v)$. Under the bistable condition $a_1>1$ and $a_2>1$, we observe that $H(u,v)=0$ is a hyperbola with one branch through $(1,0)$ and $(0,1)$ and the other branch through $(0,0)$. 
Let
\begin{eqnarray}
\mathbf{Q}_{\lambda} & = & \Big\{ (u,v) \;\Big|\; \alpha\,u+d\,\beta\,v \le \lambda,\; u,v\ge 0\Big\};\\
\mathbf{R} & = & \Big\{ (u,v) \;\Big|\; H(u,v)\ge0,\; u,v\ge 0\Big\}.
\end{eqnarray}
It is readily seen from the proof of Proposition~\ref{prop: lower bed} that a stronger lower bound can be found if we determine $\lambda_2$ in the first step for the construction of the N-barrier by 
\begin{equation}\label{eqn: representation of lambda2}
\lambda_2=\sup_{\mathbf{Q}_\lambda \subset \mathbf{R}} \lambda.
\end{equation}
To determine $\lambda_2$ given by \eqref{eqn: representation of lambda2}, we find the tangent line with the slope $-\frac{\alpha}{d\,\beta}$ to $\mathcal{L}$ at a given point on $\mathcal{L}$. To this end, we first 
solve $v=v(u)$ from $H(u,v)=0$ to get
\begin{equation}\label{eqn: v(u)}
v(u)=\frac{-(\alpha a_1u+\beta k(a_2u-1))+\sqrt{(\alpha a_1u+\beta k(a_2u-1))^2-4\alpha\beta ku(u-1)}}{2\beta k}.
\end{equation}
A straightforward calculation yields
\begin{equation}\label{eqn: dv/du}
\frac{dv(u)}{du}=\frac{-(\alpha a_1+\beta ka_2)+\frac{(\alpha a_1u+\beta k(a_2u-1))(\alpha a_1+\beta ka_2)-2\alpha\beta k(2u-1)}{\sqrt{(\alpha a_1u+\beta k(a_2u-1))^2-4\alpha\beta ku(u-1)}}}{2\beta k}.
\end{equation}
It immediately follows that the slope of the tangent line to $H(u,v)=0$ at the point $(0,1)$ (respectively, $(1,0)$) is $\frac{dv(0)}{du}=\frac{-\alpha(a_1-1)-\beta ka_2}{\beta k}$ (respectively, $\frac{dv(1)}{du}=\frac{-\alpha}{\alpha a_1+\beta k(a_2-1)}$). Under the bistable condition $a_1>1$ and $a_2>1$, we easily verify $\frac{-\alpha(a_1-1)-\beta ka_2}{\beta k}<\frac{-\alpha}{\alpha a_1+\beta k(a_2-1)}$. Noting that the slope of the line $\alpha\,u+d\,\beta\,v=\lambda_2$ is $-\frac{\alpha}{d\,\beta}$, we are led to the following three cases:
\begin{itemize}
\item[(i)] When $\frac{-\alpha}{d\beta}\le\frac{-\alpha(a_1-1)-\beta ka_2}{\beta k}$, we determine the line $\alpha\,u+d\,\beta\,v=\lambda_2$ so that it passes through $(0,1)$ and hence $\lambda_2=d\,\beta$. Note that in this case $z_2$ may be $+\infty$ in the proof of Proposition~\ref{prop: lower bed in NBMP for 2 species} (Figure~\ref{fig: NB tangent line} (i)).

\item[(ii)] When $\frac{-\alpha}{d\beta}\ge\frac{-\alpha}{\alpha a_1+\beta k(a_2-1)}$, we determine the line $\alpha\,u+d\,\beta\,v=\lambda_2$ so that it passes through $(1,0)$ and hence $\lambda_2=\alpha$. Note that in this case $z_1$ may be $-\infty$ in the proof of Proposition~\ref{prop: lower bed in NBMP for 2 species} (Figure~\ref{fig: NB tangent line} (ii)).

\item[(iii)] When $\frac{-\alpha(a_1-1)-\beta ka_2}{\beta k}<\frac{-\alpha}{d\beta}<\frac{-\alpha}{\alpha a_1+\beta k(a_2-1)}$, we determine the line $\alpha\,u+d\,\beta\,v=\lambda_2$ so that it is tangent to the curve $v=v(u)$ at some point in the first quadrant of the $uv$-plane. By \eqref{eqn: dv/du}, we have
\begin{equation}\label{eqn: dv/du case iii}
\frac{-(\alpha a_1+\beta ka_2)+\frac{(\alpha a_1u+\beta k(a_2u-1))(\alpha a_1+\beta ka_2)-2\alpha\beta k(2u-1)}{\sqrt{(\alpha a_1u+\beta k(a_2u-1))^2-4\alpha\beta ku(u-1)}}}{2\beta k}=\frac{-\alpha}{d\beta}
\end{equation}
or 
\begin{equation}\label{eqn:solve_u}
\frac{Au^2+Bu+C}{Du^2+Eu+J}=G,
\end{equation}
where
\begin{alignat}{3}\nonumber
A &=\left(X^2-4\alpha\beta k\right)^2, & \quad B &=(X^2-4\alpha\beta k)(-\beta kX+2\alpha\beta k),
\\\notag
C &=\left(-\beta kX+2\alpha\beta k\right)^2, & \quad D &=X^2-4\alpha\beta k, 
\\\notag
E &=-2\beta kX+4\alpha\beta k, & \quad J &=\beta^2 k^2,
\\\notag   
G &=\left(X-2\alpha k\,d^{-1}\right)^2, & \quad X &=\alpha a_1+\beta ka_2.
\end{alignat}
Then $u$ is solved from \eqref{eqn:solve_u} by
\begin{equation}\label{eqn: expression of u}
u=\frac{-(B-EG)\pm\sqrt{(B-EG)^2-4(A-DG)(C-JG)}}{2(A-DG)}.
\end{equation}
Consequently, $\lambda_2=\alpha\,u+d\,\beta\,v(u)$ is determined by \eqref{eqn: v(u)} and \eqref{eqn: expression of u}, where $u$ is given by \eqref{eqn: expression of u} with u satisfying $0<u<1$ and \eqref{eqn: dv/du case iii} (Figure~\ref{fig: NB tangent line} (iii)).


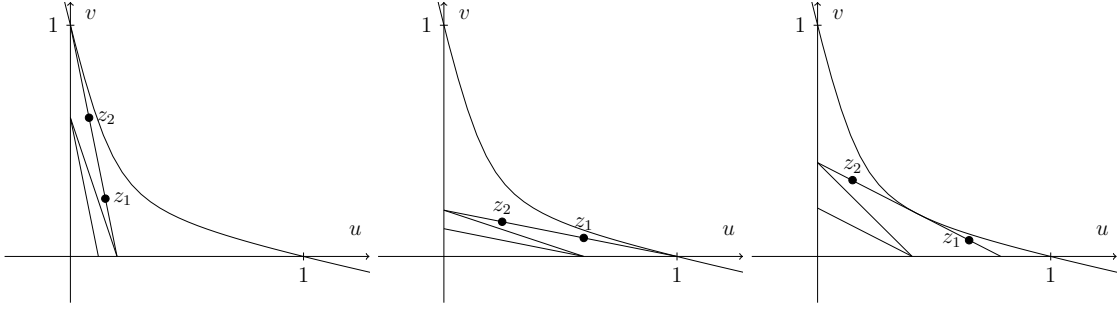
\begin{figure}
\begin{minipage}[t]{0.3\linewidth}
\centering
\begin{tikzpicture}[scale=0.70]
\begin{axis}[axis equal,
xmin=-0.1, xmax=1.1, ymin=-0.2, ymax=1.1,
 axis lines=middle,
 axis line style={->},
 tick style={color=black},
 xtick=1,
 ytick=1,
 xlabel=$u$,x label style={at={(1,0.2)}},
 ylabel=$v$,y label style={at={(0.2,1)}},
]
  \addplot[samples=100, domain=-2:2] 
	({x}, {1/2*(-5*x+1+(21*x^2-6*x+1)^0.5)});
\addplot[samples=100, domain=0:1/5] 
	({x}, {-5*x+1});
\addplot[samples=100, domain=0:1/5] 
	({x}, {-3*(x-1/5)});
\addplot[samples=100, domain=0:3/25] 
	({x}, {-5*x+3/5});
  \addplot[color=black,mark=*] coordinates {(0.08,0.6)};
\addplot[color=black,mark=*] coordinates {(0.15,0.25)};
 \node[right] at (axis cs:0.08,0.6) {$z_2$};
  \node[right] at (axis cs:0.15,0.25) {$z_1$};
\end{axis}
\end{tikzpicture}
\end{minipage}
\begin{minipage}[t]{0.3\linewidth}
\centering
\begin{tikzpicture}[scale=0.70]
\begin{axis}[axis equal,
xmin=-0.1, xmax=1.1, ymin=-0.2, ymax=1.1,
 axis lines=middle,
 axis line style={->},
 tick style={color=black},
 xtick=1, 
 ytick=1, 
 xlabel=$u$,x label style={at={(1,0.2)}},
 ylabel=$v$,y label style={at={(0.2,1)}},
]
  \addplot[samples=100, domain=-2:2] 
	({x}, {1/2*(-5*x+1+(21*x^2-6*x+1)^0.5)});
\addplot[samples=100, domain=0:1] 
	({x}, {-1/5*x+1/5});
\addplot[samples=100, domain=0:3/5] 
	({x}, {-1/3*x+1/5)});
\addplot[samples=100, domain=0:3/5] 
	({x}, {-1/5*(x-3/5)});
  \addplot[color=black,mark=*] coordinates {(0.6,0.08)};
\addplot[color=black,mark=*] coordinates {(0.25,0.15)};
 \node[above] at (axis cs:0.6,0.08) {$z_1$};
  \node[above] at (axis cs:0.25,0.15) {$z_2$};
\end{axis}
\end{tikzpicture}
\end{minipage}
\begin{minipage}[t]{0.3\linewidth}
\centering
\begin{tikzpicture}[scale=0.70]
\begin{axis}[axis equal,
xmin=-0.1, xmax=1.1, ymin=-0.2, ymax=1.1,
 axis lines=middle,
 axis line style={->},
 tick style={color=black},
 xtick=1, 
 ytick=1, 
 xlabel=$u$,x label style={at={(1,0.2)}},
 ylabel=$v$,y label style={at={(0.2,1)}},
]
  \addplot[samples=100, domain=-2:2] 
	({x}, {1/2*(-5*x+1+(21*x^2-6*x+1)^0.5)});
\addplot[samples=100, domain=0:(4/7-2/7^0.5)/(-5/2+3*7^0.5/4)+3/7] 
	({x}, {(-5/2+3*7^0.5/4)*(x-3/7)-(4/7-2/7^0.5)});
\addplot[samples=100, domain=0:(-5/2+3*7^0.5/4)*(-3/7)-(4/7-2/7^0.5)] 
	({x}, {-(x-((-5/2+3*7^0.5/4)*(-3/7)-(4/7-2/7^0.5)))});
\addplot[samples=100, domain=0:(-5/2+3*7^0.5/4)*(-3/7)-(4/7-2/7^0.5)] 
	({x}, {(-5/2+3*7^0.5/4)*(x-((-5/2+3*7^0.5/4)*(-3/7)-(4/7-2/7^0.5)))});
  \addplot[color=black,mark=*] coordinates {(0.65,0.07)};
\addplot[color=black,mark=*] coordinates {(0.15,0.33)};
 \node[left] at (axis cs:0.65,0.07) {$z_1$};
  \node[above] at (axis cs:0.15,0.33) {$z_2$};
\end{axis}
\end{tikzpicture}
\end{minipage}
\caption{N-barrier for cases (i), (ii) and (iii) (from the left to the right)}
\label{fig: NB tangent line}
\end{figure}

\end{itemize}

\vspace{2mm}
\setcounter{equation}{0}
\setcounter{figure}{0}
\section{Nonexistence of four species waves}\label{sec: nonexistence}
\vspace{2mm}


With the aid of Corollary~\ref{cor: NBMP for 3 species}, we establish a nonexistence result for \textbf{(LV4)} in this section. Recall \textbf{(LV4)} is as follows:
\begin{equation}\nonumber
\textbf{(LV4)}
\begin{cases}
\vspace{3mm}
d_i(u_i)_{xx}+\theta(u_i)_{x}+u_i(\sigma_i-c_{i1}\,u_1-c_{i2}\,u_2-c_{i3}\,u_3-c_{i4}\,u_4)=0,\; x\in\mathbb{R},\; i=1,...,4, \\
(u_1,u_2,u_3,u_4)(-\infty)=(\frac{\sigma_1}{c_{11}},0,0,0),\quad 
(u_1,u_2,u_3,u_4)(\infty)=(0,\frac{\sigma_2}{c_{22}},0,0).
\end{cases}
\end{equation}

\begin{thm}[\textbf{Nonexistence of four species waves}]\label{thm: Nonexistence 4 species}
Assume that the following hypotheses hold:
\begin{itemize}
\item [$\mathbf{[H1]}$] $\tilde{\sigma}_i:=\sigma_i-c_{i4}\,\sigma_4\,c_{44}^{-1}>0$ for $i=1,2,3$;
\item [$\mathbf{[H2]}$] 
$\displaystyle\min_{i=1,2,3}
\Big(
\alpha^{\ast}_i\displaystyle\min_{j=1,2,3}\frac{\tilde{\sigma}_j}{c_{ji}}
\Big)\,
\Big(
\displaystyle\min_{i=1,2,...,n} d_i 
\Big)
\Big(
\displaystyle\max_{i=1,2,...,n} d_i
\Big)^{-1}\geq\sigma_4$, where $\alpha^{\ast}_i=c_{4i}$, $i=1,2,3$. 
\end{itemize}
Then \textbf{(LV4)} has no positive solution $(u_1(x),u_2(x),u_3(x),u_4(x))$.
\end{thm}

\begin{proof}
We prove by contradiction. Suppose to the contrary that there exists a solution $(u_1(x),u_2(x),u_3(x),u_4(x))$ to \textbf{(LV4)}. It follows from the fact $u_4(x)>0$ for $x\in\mathbb{R}$ and $u_4(\pm\infty)=0$, that there exists $x_0\in\mathbb{R}$ such that $\max_{x\in\mathbb{R}} u_4(x)=u_4(x_0)>0$, $u_4''(x_0)\le 0$, and $u_4'(x_0)=0$. Due to $d_4\,(u_4)_{xx}+\theta \,(u_4)_x+u_4\,(\sigma_4-c_{41}\,u_1-c_{42}\,u_2-c_{43}\,u_3-c_{44}\,u_4)=0$, we obtain
\begin{equation}\label{eqn: w(???) ???>0}
\sigma_4-c_{41}\,u_1(x_0)-c_{42}\,u_2(x_0)-c_{43}\,u_3(x_0)-c_{44}\,u_4(x_0)\ge 0,
\end{equation}
which gives
\begin{equation}
u_4(x)\leq u_4(x_0)\le \frac{1}{c_{44}}\big(\sigma_4-c_{41}\,u_1(x_0)-c_{42}\,u_2(x_0)-c_{43}\,u_3(x_0)\big)<\frac{\sigma_4}{c_{44}},\;x\in\mathbb{R}.
\end{equation}
As a result, we have
\begin{equation}\label{eqn: nonexistence diff ineq <0}
\begin{cases}
\vspace{3mm}
d_1\,(u_1)_{xx}+\theta \,(u_1)_x+u_1\,(\sigma_1-c_{14}\,\sigma_4\,c_{44}^{-1}-c_{11}\,u_1-c_{12}\,u_2-c_{13}\,u_3)\leq0, \quad x\in\mathbb{R}, \\
\vspace{3mm}
d_2\,(u_2)_{xx}+\theta \,(u_2)_x+u_2\,(\sigma_2-c_{24}\,\sigma_4\,c_{44}^{-1}-c_{21}\,u_1-c_{22}\,u_2-c_{23}\,u_3)\leq0, \quad x\in\mathbb{R}, \\
d_3\,(u_3)_{xx}+\theta \,(u_3)_x+u_3\,(\sigma_3-c_{34}\,\sigma_4\,c_{44}^{-1}-c_{31}\,u_1-c_{32}\,u_2-c_{33}\,u_3)\leq0, \quad x\in\mathbb{R}.
\end{cases}
\end{equation}
Because of $\mathbf{[H1]}$, we apply Corollary~\ref{cor: NBMP for 3 species} to the last three inequalities, and obtain a lower bound of $c_{41}\,u_1(x)+c_{42}\,u_2(x)+c_{43}\,u_3(x)$:
\begin{equation}
c_{41}\,u_1(x)+c_{42}\,u_2(x)+c_{43}\,u_3(x)\geq \min_{i=1,2,3}
\Big(
\alpha^{\ast}_i\min_{j=1,2,3}\frac{\tilde{\sigma}_j}{c_{ji}}
\Big)\,
\Big(
\min_{i=1,2,...,n} d_i 
\Big)
\Big(
\max_{i=1,2,...,n} d_i
\Big)^{-1},\;x\in\mathbb{R}.
\end{equation}
The hypothesis $\mathbf{[H2]}$ then yields
\begin{equation}
c_{41}\,u_1(x)+c_{42}\,u_2(x)+c_{43}\,u_3(x)\geq\sigma_4,\;x\in\mathbb{R},
\end{equation}
which contradicts \eqref{eqn: w(???) ???>0}. This completes the proof.

\end{proof}

\textbf{Biological interpretation of Theorem~\ref{thm: Nonexistence 4 species}}:
When other parameters are fixed, it is easy for $\mathbf{[H1]}$ and $\mathbf{[H2]}$ to hold true as long as $\sigma_4$ is sufficiently small. Biologically, this means that when the intrinsic growth rate $\sigma_4$ of $u_4$ is sufficiently small, the four species $u_1$, $u_2$, $u_3$ and $u_4$ cannot coexist in the ecological system \textbf{(LV4)} under certain parameter regimes.



\vspace{2mm}

\textbf{Acknowledgments.} The authors wish to express sincere gratitude to
 Dr. Tom Mollee for his careful reading of the manuscript and valuable suggestions and comments
to improve the readability and accuracy of the paper. The research of L.-C. Hung is partly supported by the grant 104EFA0101550 of Ministry of Science and Technology, Taiwan. 




\vspace{2mm}











\end{document}